\documentclass[a4paper,12pt,reqno]{amsart}
\usepackage{amssymb,amsthm ,amsmath,mathtools}
\usepackage{ifthen}
\usepackage{graphicx}
\usepackage{float}
\usepackage{color} 
\usepackage{xcolor}
\usepackage{diagbox}
\usepackage{caption}
\usepackage{subcaption}
\usepackage{textcomp,multirow}
\theoremstyle{plain}
\newtheorem{thm}{Theorem}[section]
\newtheorem*{thm*}{Theorem}
\newtheorem{conj}{Conjecture}[section]

\numberwithin{equation}{section}

\newtheorem{lem}{Lemma}[section]

\newtheorem{rem}{Remark}[section]

\theoremstyle{definition}
\newcounter {own}
\def\theown {\thesection  .\arabic{own}}

\newenvironment{pf}[1][]{%
 \vskip 3mm
 \noindent
 \ifthenelse{\equal{#1}{}}%
  {{\slshape Proof. }}%
  {{\slshape #1.} }%
 }%
{\qed\bigskip}
\graphicspath{ {./images/} }
\DeclarePairedDelimiter\ceil{\lceil}{\rceil}
\DeclarePairedDelimiter\floor{\lfloor}{\rfloor}

\DeclareMathOperator*{\supp}{\supp}
\DeclareMathOperator*{\rank}{rank}

\newcounter{alphabet}

\newcommand{\ds}{\displaystyle}

\newcounter{minutes}
\divide\time by 60
\newcounter{hours}
\multiply\time by 60 \addtocounter{minutes}{-\time}

\begin{document}
\bibliographystyle{amsplain}
\title{Obstructions for Gabor frames of the second order B-spline}

\thanks{
File:~AntonyRiya6.tex,
          printed: 2023-01-10,
          \thehours.\ifnum\theminutes<10{0}\fi\theminutes}

\author{Riya Ghosh}
\author{A. Antony Selvan$^\dagger$}

\address{Riya Ghosh, Indian Institute of Technology Dhanbad, Dhanbad-826 004, India.}
\email{riya74012@gmail.com}
\address{A. Antony Selvan, Indian Institute of Technology Dhanbad, Dhanbad-826 004, India.}
\email{antonyaans@gmail.com}

\keywords{B-splines, Gabor frames, Laurent operator, shift-invariant spaces, Zak transform, Zibulski–Zeevi matrix.\\
$^\dagger$ {\tt Corresponding author}
}
\maketitle
\pagestyle{myheadings}
\markboth{Riya Ghosh and A. Antony Selvan}{Obstructions for Gabor frames of the second order B-spline}
\begin{abstract}
For a window 
 $g\in L^2(\mathbb{R})$, the subset of all lattice parameters
$(a, b)\in \mathbb{R}^2_+$ such that $\mathcal{G}(g,a,b)=\{e^{2\pi ib m\cdot}g(\cdot-a k) : k, m\in\mathbb{Z}\}$ forms a frame for $L^2(\mathbb{R})$ is known as the frame set of $g$. In time-frequency analysis, determining the Gabor frame set for a given window is a challenging open problem. In particular, the frame set for B-splines has many obstructions.
Lemvig and Nielsen in \cite{counter} conjectured that if 
\begin{align} 
   a_0=\dfrac{1}{2m+1},~ b_0=\dfrac{2k+1}{2},~k,m\in \mathbb{N},~k>m,~a_0b_0<1,\nonumber
\end{align}
then the Gabor system $\mathcal{G}(Q_2, a, b)$ of the second order B-spline $Q_2$ is not a frame along the hyperbolas
    \begin{align}
    ab=\dfrac{2k+1}{2(2m+1)},\text{ for }b\in \left[b_0-a_0\dfrac{k-m}{2}, b_0+a_0\dfrac{k-m}{2}\right],\nonumber
\end{align}
for every $a_0$, $b_0$. 
Nielsen  in \cite {Nielsenthesis} also  conjectured that $\mathcal{G}(Q_2, a,b)$ is not a frame for 
$$a=\dfrac{1}{2m},~b=\dfrac{2k+1}{2},~k,m\in \mathbb{N},~k>m,~ab<1\text{ with }\gcd(4m,2k+1)=1.$$ In this paper, we prove that both conjectures are true.
\end{abstract}
\section{Introduction}
One of the fundamental problems in Gabor analysis is to determine the values of $a, b>0$ such that the Gabor system  $\mathcal{G}(g, a, b)=\{e^{2\pi ib m\cdot}g(\cdot-a k): k, m\in\mathbb{Z}\}$ for a nonzero window function $g\in L^2(\mathbb{R})$ forms a frame for $L^2(\mathbb{R})$, \textit{i.e.}, there exist two positive constants $A$, $B$ such that 
\begin{equation}\label{Gaborframe}
A\|f\|^2\leq\displaystyle\sum_{n\in\mathbb{Z}}|\langle f,e^{2\pi ib m\cdot}g(\cdot-a k)\rangle|^2\leq B\|f\|^2,
\end{equation}
for every $f\in L^2(\mathbb{R})$. The constants $A$ and $B$ are called frame bounds.  For the window function $g$, the set $\mathcal{F}(g)= \left\{(a, b) \in\mathbb{R}^2_+ : \mathcal{G}(g, a,b)~ \text{is a frame}\right\} $ is called the frame set of $g$. The fundamental density theorem asserts that 
$$\mathcal{F}(g)\subseteq \{(a, b) \in\mathbb{R}^2_+: ab\le 1\}$$ 
(see \cite{hedtg}). In addition, if $g$ is in the modulation space $M^1(\mathbb{R})$, Feichtinger and Kaiblinger \cite{HGF1} proved that 
$\mathcal{F}(g)$ is an open subset of $\mathbb{R}^2_+$ and Balian-Low theorem states that $\mathcal{F}(g)\subseteq \{(a, b) \in\mathbb{R}^2_+: ab< 1\}$ \cite{balian, rdbalian}. The frame set is completely characterized only for a few windows: the Gaussian, the hyperbolic secant, the two-sided exponential, the one-sided exponential, the characteristic function, the totally positive functions of finite type $\ge 2$ or of Gaussian type, and the Herglotz functions (see \cite{gfrf, abc, duke, stsis, totally, when, whfb, hsyg, tgfcd, fblyu, dtsibf2}).

Gabor frames with compactly supported windows play a significant role in Gabor analysis due to their inherent time-frequency localization and they have received a lot of attention in recent years, 
see \cite{gwole, sign, mystery} and references therein. For a function 
$g$ with $\text{supp }g\subseteq [0, L]$, 
$\mathcal{F}(g)$ is always a subset of $\left\{(a, b)\in\mathbb{R}_{+}^2: ab \le 1~\text{and}~a \le L\right\}.$ 
The complete frame set of the first order B-spline \cite{abc, when} and the Haar function \cite{haar} are extremely complicated. The frame set is much less clear for the higher order B-splines.
The B-splines are defined inductively as follows:
\begin{eqnarray}\label{Qmformula}
Q_1(x):=\chi_{\left[-\tfrac{1}{2},\tfrac{1}{2}\right]}(x)~\text{and}~Q_{n+1}(x):=(Q_n*Q_1)(x),~n\geq1.
\end{eqnarray}
 Since $Q_n\in M^1(\mathbb{R})$ for $n \ge 2$, $\mathcal{F}(Q_2)$ is an open subset of $\mathbb{R}^2$ but a full characterization of this set still an open problem. In recent years, some progress in \cite{fscf, b3spline, ofsb, sign, olec, ourgabor, counter} has been made to characterize the frame set for B-splines. In \cite{mystery},
Gr\"{o}chenig initially conjectured that 
$$\mathcal{F}(Q_n) =\{(a, b)\in \mathbb{R}^2_{+} : ab < 1, a < n, b\ne 2, 3,\dots\},$$
but Lemvig and Nielsen \cite{counter} disproved it. They obtained a family of obstructions for the Gabor frame property of B-splines stating that $\mathcal{G}(Q_n, a, b)$ is not a frame for $L^2(\mathbb{R})$ when $b>3/2$, $ab=p/q$ with $\gcd(p,q)=1$ and $\big|b-\floor{b+0.5}\big|\le{1}/{nq}$. They proved that if $ab = 5/6$ and $b \in[\tfrac{7}{3} , \tfrac{8}{3}]$, then $\mathcal{G}(Q_2, a, b)$ is not a frame for $L^2(\mathbb{R})$. Further, they considered more obstructions for the Gabor frame property of the second order B-spline and conjectured the following.\
\begin{conj}[\textbf{Lemvig and Nielsen}]\cite{counter}\label{conjecture}
Let 
\begin{align}\label{a0b0} 
   a_0=\dfrac{1}{2m+1},~ b_0=\dfrac{2k+1}{2},~k,m\in \mathbb{N},~k>m,~a_0b_0<1.
\end{align}
Then the Gabor system $\mathcal{G}(Q_2, a,b)$ is not a frame along the hyperbolas
    \begin{align}\label{ab}
    ab=\dfrac{2k+1}{2(2m+1)},\text{ for }b\in \left[b_0-a_0\dfrac{k-m}{2}, b_0+a_0\dfrac{k-m}{2}\right],
\end{align}
for every $a_0$, $b_0$.
\end{conj}

\begin{conj}[\textbf{Nielsen}]\cite{Nielsenthesis}\label{conjecture2}
The Gabor system $\mathcal{G}(Q_2, a,b)$ is not a frame for 
$$a=\dfrac{1}{2m},~b=\dfrac{2k+1}{2},~k,m\in \mathbb{N},~k>m,~ab<1\text{ with }\gcd(4m,2k+1)=1.$$
\end{conj}
Gr\"{o}chenig \cite{partition} verified that $\mathcal{G}(Q_2,a_0,b_0)$ is not a frame for $L^2(\mathbb{R})$. In this paper, we prove the following.
\begin{thm}\label{main thm}
Conjectures \ref{conjecture} and \ref{conjecture2} are true.
\end{thm}
 
Nielsen mentioned in her thesis \cite{Nielsenthesis} that there appeared to be an interval along the hyperbolas corresponding to the points in Conjecture \ref{conjecture2}, where the Gabor systems $\mathcal{G}(Q_2, a, b)$ are not frames. However, she was not able to determine the range of $b$. On the basis of numerical experiments, we expect that if
$\alpha_0=\dfrac{1}{2m},~ \beta_0=\dfrac{2k+1}{2},~k,m\in \mathbb{N},~k>m,~\alpha_0\beta_0<1\text{ with }\gcd{(4m,2k+1)}=1$, then $\mathcal{G}(Q_2, a,b)$ is not a frame along the hyperbolas
\begin{align}\label{abnew}
    ab=\dfrac{2k+1}{4m},\text{ for }b\in \left[\beta_0-\alpha_0\dfrac{k-m}{2}, \beta_0+\alpha_0\dfrac{k-m}{2}\right],
 \end{align}
for every $\alpha_0$, $\beta_0$. But we are not able to prove it.

The paper is organized as follows. In Section 2, we discuss characterizations for the Gabor frames with rational density in the literature and provide a new characterization using the Laurent operator technique. In Section 3, we prove some important Lemmas to prove our main result. Finally, the proof of Theorem \ref{main thm} is included in Section 5. 

\section{Characterizations of Gabor frames with rational density}
Let us introduce the shift-invariant space of a generator $g$ defined by 
$$V_h(g):=\left\{f\in L^2(\mathbb{R}): f(\cdot)=\sum\limits_{k\in\mathbb{Z}}d_k g(\cdot-hk)~\text{for some}~ (d_k)\in \ell^2(\mathbb{Z})\right\},~h>0.$$
Recall that $g$ is said to be a stable generator for $V_h(g)$ if $\{g(\cdot-hk) : k \in \mathbb{Z}\}$ is a Riesz basis for $V_h(g)$, \textit{i.e.},
$\overline{span}\{g(\cdot-hk): k\in\mathbb{Z}\}=V_h(g)$ and there exist constants $A$, $B>0$ such that
\begin{equation}\label{rieszbasis}
A\sum_{k\in\mathbb{Z}}|d_k|^2\leq\big\|\sum_{k\in\mathbb{Z}}d_k g(\cdot-hk)\big\|^2\leq
B\sum_{k\in\mathbb{Z}}|d_k|^2,
\end{equation}
for all $(d_k)\in\ell^2(\mathbb{Z})$.    
A set $\Lambda=\{x_{n}:n\in\mathbb{Z}\}$ of real numbers is said to be a set of stable sampling for $V_h(g)$ if there exist constants $A$, $B>0$ such that
\begin{eqnarray}\label{pap3eqn2.7}
A\| f\|^2\leq\ds\sum\limits_{n\in\mathbb{Z}}\ds|f(x_n)|^2 \leq B\| f\|^2, 
\end{eqnarray}
for all $f\in V_h(g)$. The numbers $A$ and $B$ are called sampling bounds. The following theorem establishes a fundamental link between Gabor analysis and the theory of sampling in shift-invariant spaces.
\begin{thm}\label{main}\cite{J95, Ron}.
Let $g$ be a stable generator for $V_{1/b}(g)$. Then
$\mathcal{G}(g,a, b)$ forms a frame for $L^2(\mathbb{R})$ if and only if $x+a\mathbb{Z}$ is a set of stable sampling for $V_{1/b}(g)$ with uniform constants independent of $x\in \mathbb{R}$.
\end{thm}

Given a window function $f$ in the Schwartz class $S(\mathbb{R})$, consider the short-time Fourier transform  of $g$ with respect to $f$:
$$V_{f} g(x,w) =\int_{-\infty}^\infty
g(t)f (x - t)e^{2\pi iwt} dt.$$
A function $g$ belongs to the modulation space $M^1(\mathbb{R})$ if
     $$\int_{-\infty}^\infty\int_{-\infty}^\infty \big|V_f g(x,w)\big|~ dx dw < \infty$$
     for some non-trivial function $f \in S(\mathbb{R})$.
The Zak transform is a major tool in the analysis of Gabor frames with rational density.
For a fixed $\alpha >0$, it is defined as
$$(\mathcal{Z}_{\alpha}g)(x,t)=\sum\limits_{n\in\mathbb{Z}}g\left(x-\alpha n\right)e^{2\pi i\alpha nt},~g\in L^2(\mathbb{R}).$$
If $g\in M^1(\mathbb{R})$, then $\mathcal{Z}_{\alpha} g$ is continuous.

Given $p, q \in\mathbb{N}$ with $p\le q$, $\text{gcd}(p,q)=1$, and $ab=\tfrac{p}{q}$, we define a $p\times q$ matrix 
\begin{align}\label{zeevi}
    \Psi_g(x, t)=\left[(\mathcal{Z}_{a}g)\left(x+\dfrac{ak}{p},t-bl\right)e^{2\pi ikl/q}\right]_{k=0,~l=0}^{p-1,~q-1}.
\end{align}
The matrix $\Psi_g(x, t)$ is called the Zibulski– Zeevi matrix.
Let us define another $p\times q$ matrix 
\begin{align}\label{counterexample}
     \Theta_g(x, t)=
\left[\sum\limits_{n\in \mathbb{Z}}
g(x + aqn + al+ k/b) e^{-2\pi iaqnt}\right]_{k=0,~l=0}^{p-1,~q-1}.
\end{align}
The matrix $\Theta_g(x,t)$ is a variant of the
Zibulski–Zeevi matrix and appears in \cite{cocompact, Lyu}. 

The following theorem characterizes Gabor frames with rational density.  
\begin{thm}\cite{Lyu, zeevi}.
Assume that $ab=\tfrac{p}{q}$ with $\gcd(p,q)=1$ and $g\in M^1
(\mathbb{R})$. Then the following statements are equivalent:
\begin{enumerate}
    \item [$(i)$] $\mathcal{G}(g, a, b)$ is a frame for $L^2(\mathbb{R}).$
    \item[$(ii)$] $\Psi_g(x, t)$ has rank $p$ for all $x, t\in \mathbb{R}$.
    \item[$(iii)$] $\Theta_g(x, t)$ has rank $p$ for all $x, t\in \mathbb{R}$.
\end{enumerate}
\end{thm}

We now give another
characterization for Gabor frames with rational density using block Laurent operators.

Let $\mathcal{A}$ be a bounded linear operator on $\ell^2(\mathbb{Z})$ with
the associated matrix given by $\mathcal{A}=[a_{rs}]$. Then $\mathcal{A}$ is said to be a \textit{Laurent operator} if $a_{r-k,s-k}=a_{r,s}$, for every $r$, $s$, $k$
$\in\mathbb{Z}$. For $m\in L^\infty(\mathbb{T})$, let us define $\mathcal{M}:L^2(\mathbb{T})\to
L^2(\mathbb{T})$ by \[(\mathcal{M}f)(x):=m(x)f(x), ~\text{for}~a.e.~x\in\mathbb{T},\]
where $\mathbb{T}$ is the circle group identified with the unit interval  $[0,1)$. Then the operator $\mathcal{A}=\mathcal{F}\mathcal{M}\mathcal{F}^{-1}$ is a Laurent operator defined by the matrix $[\widehat{m}(r-s)]$, where $\mathcal{F}:L^2(\mathbb{T})\to \ell^2(\mathbb{Z})$ is the Fourier transform defined by $\mathcal{F}f=\widehat{f},~\widehat{f}(n)=\int_{0}^{1} f(x) e^{-2\pi in x}dx,~n\in\mathbb{Z}$. In this case, we say that $\mathcal{A}$ is a Laurent operator defined by the symbol $m$. If it is invertible, then $\mathcal{A}^{-1}$ is also a Laurent operator with the symbol $\tfrac{1}{m}$ and the matrix of $\mathcal{A}^{-1}$ is given by
$$\mathcal{A}^{-1}=\left[\widehat{\left(\tfrac{1}{m}\right)}(r-s)\right].$$

Let $\mathcal{L}$ be a bounded linear operator on $\ell^2(\mathbb{Z})^m$ 
which is represented by an $m\times m$ matrix whose entries are doubly infinite matrices, \textit{i.e.},
\begin{eqnarray*}
 \mathcal{L}=\left(
\begin{array}{ccccccc}
 L_{11} & L_{12}  & \cdots & L_{1m} \\
L_{21} & L_{22}  &  \cdots & L_{2m} \\
\vdots & \vdots  & \ddots & \vdots\\
L_{m1} & L_{m2}  & \cdots & L_{mm}  
 \end{array}
\right), ~ L_{rs} ~\text{is a doubly infinite matrix}.
\end{eqnarray*}
Then $\mathcal{L}$ is said to be a \textit{block Laurent operator} if each $L_{ij}$'s are Laurent operators on $\ell^2(\mathbb{Z})$. 

For an $m\times m$ matrix-valued function $\Phi \in L_{m\times m}^1(\mathbb{T})$, the Fourier coefficient $\widehat{\Phi}(k)$ of $\Phi$ is an $m\times m$  matrix 
$$ \widehat \Phi (k):= \int\limits_{\mathbb{T}}  e^{-2\pi i k  x}\Phi(x) dx, ~ k \in \mathbb{Z},$$
whose $(i,j)$-th entry is equal to the Fourier coefficient of the $(i,j)$-th entry of $\Phi$. It is well known that for every block Laurent operator $\mathcal{L}$, there exists a function $\Phi \in L^{\infty}_{m\times m}(\mathbb{T})$ such that $\mathcal{L}=[\widehat{\Phi}(i-j)]$ (see \cite[p. 565-566]{GoGoKa1}). In this case, we say that $\mathcal{L}$ is a block Laurent operator defined by $\Phi$. Since the class of Laurent operators is closed under addition, multiplication, and multiplication is commutative,
$$\det \mathcal{L}:=\sum\limits_{\tau} (sgn~\tau)L_{1\tau_1}\cdots L_{m{\tau}_m},$$
is a well-defined Laurent operator  on $\ell^2(\mathbb{Z})$ with the symbol $\det \Phi$. If $\mathcal{L}$ is a block Laurent operator with the symbol $\Phi\in L^\infty_{m\times m}(\mathbb{T}),$ then $\mathcal{L}^*$ is a block Laurent operator with the symbol $\Phi^*,$ where $\Phi^*$ denotes the adjoint of $\Phi.$ Consequently $\mathcal{L}^{*}\mathcal{L}$ is a block Laurent operator with the symbol $\Phi^*\Phi.$ 
\begin{thm}\label{IB} \cite{GoGoKa1}.
Let $\mathcal{L}=[L_{rs}]=[\widehat{\Phi}(i-j)]$ denote the matrix of a function $\Phi\in L^2_{m\times m}(\mathbb{T})$. Then $\mathcal{L}$ is an invertible block Laurent operator with the symbol $\Phi\in L^\infty_{m\times m}(\mathbb{T})$ if and only if there exist  two positive constants $A$ and $B$ such that $$A\leq |\det\Phi(x)|\leq B,~\text{for}~a.e. ~x\in\mathbb{T}.$$
\end{thm}

Let $ab=\tfrac{p}{q}\in \mathbb{Q}$ with $\gcd(p,q)=1$ and $g\in M^1(\mathbb{R})$. For $f\in V_{\tfrac{1}{b}}(g)$, we have
\begin{eqnarray}\label{f1b}
f(t)=\sum\limits_{k\in\mathbb{Z}}c_kg(t-\tfrac{k}{b})=\sum\limits_{k\in\mathbb{Z}}\sum\limits_{n=0}^{p-1}c_{pk+n}g\left(t-aqk-\tfrac{n}{b}\right).
\end{eqnarray}
Consider the sampling set $X = x+a\mathbb{Z}=\bigcup\limits_{s=0}^{q-1}\big\{x+aql+as:~l\in\mathbb{Z}\big\},~x\in\mathbb{R}$ and construct an infinite system 
\begin{eqnarray*}
f(x+aql+as)=\sum\limits_{k\in\mathbb{Z}}\sum\limits_{n=0}^{p-1}c_{pk+n}g\left(x+aq(l-k)+as-\tfrac{n}{b}\right),
\end{eqnarray*}
where  $l\in \mathbb{Z}$ and $s=0,1,\dots,q-1.$ The above system can be written as 
\begin{eqnarray}
UC=F,
\end{eqnarray}
where 
\begin{eqnarray}
U=\begin{bmatrix}
U_{00}&U_{01}&\cdots&U_{0,p-1}\\
U_{10}&U_{11}&\cdots&U_{1,p-1}\\
\vdots&\vdots&\ddots&\vdots\\
U_{q-1,0}&U_{q-1,1}&\cdots&U_{q-1,p-1}
\end{bmatrix},~C=\begin{bmatrix}
    C_0\\C_1\\\vdots\\C_{p-1}
\end{bmatrix},\text{and}~F=\begin{bmatrix}
    F_0\\F_1\\\vdots\\F_{q-1}
\end{bmatrix},
\end{eqnarray}
with $$U_{sn}=\left[g
\left(x+aq(l-k)+as-\tfrac{n}{b}\right)\right]_{l,k\in\mathbb{Z}},~C_n=\{c_{pk+n}\}_{k\in\mathbb{Z}}^{T},$$
and
$$F_{s}=\big\{f(\alpha+aql+as)\big\}_{l\in\mathbb{Z}}^{T},$$
for $0\le n\le p-1$ and $0\le s\le q-1$. It is easy to check that $X=x+a\mathbb{Z}$ is an SS if and only if $U^{*}U$ is invertible on $\ell^2(\mathbb{Z})^p$. Notice that $U_{sn}$ is a Laurent operator with the symbol 
\begin{equation}\label{symbolentries}
 \Phi_{sn}(x,t)=\ds\sum\limits_{k\in\mathbb{Z}}g\left(x+aqk+as-\dfrac{n}{b}\right)e^{2\pi ikt}=(\mathcal{Z}_{aq}g)\left(x+as-\dfrac{n}{b},-\dfrac{t}{aq}\right).
\end{equation}
Therefore, $U^{*}U$ is a block Laurent operator on $\ell^2(\mathbb{Z})^{p}$ with the symbol $\Phi_g^{*}(x,t)\Phi_g(x,t)$, where 
\begin{eqnarray}\label{phi}
\Phi_g(x,t)=
\begin{bmatrix}
\Phi_{00}(x,t)&\Phi_{01}(x,t)&\cdots&\Phi_{0,p-1}(x,t)\\
\Phi_{10}(x,t)&\Phi_{11}(x,t)&\cdots&\Phi_{1,p-1}(x,t)\\
\vdots&\vdots&\ddots&\vdots\\
\Phi_{q-1,0}(x,t)&\Phi_{q-1,1}(x,t)&\cdots&\Phi_{q-1,p-1}(x,t)
\end{bmatrix}.
\end{eqnarray}
The matrix $\Phi_g(x,t)$ is closely related to $\Theta_g(x,t)$. In particular, $\Phi_g(0,0)=\Theta^{t}_g(0,0).$ 
It is well known that the Gram matrix $M^*M$ of a matrix $M$ is invertible if and only if $M$ has linearly independent columns. From Theorems \ref{main} and \ref{IB} we derive the following result.
\begin{thm}\label{framewithrank}
    Let $g\in M^1(\mathbb{R})$ be a stable generator for $V_{1/b}(g)$ and $ab\in \mathbb{Q}$. Then $\mathcal{G}(g, a, b)$ forms a frame for $L^2(\mathbb{R})$ if and only if $\Phi_g(x,t)$ defined in \eqref{phi} has linearly independent columns for all $x,t\in \mathbb{R}$.
\end{thm}
\section{Key Lemmas}
For convenience, we take $p=2k+1$ and $q=2(2m+1)$ with $\gcd(p,q)=1$. If $\gcd(p,q)\ne 1$, then we can rewrite $ab=\tfrac{2k+1}{2(2m+1)}=\tfrac{2k_1+1}{2(2m_1+1)}$, for some $k_1>m_1$ such that $\gcd(2k_1+1,2(2m_1+1))=1$.
The condition \eqref{a0b0} implies that $m+1\le k\le 2m$ and the condition \eqref{ab} implies that
$$\dfrac{(2k+1)(2m+1)+m-k}{2(2m+1)}\le b\le \dfrac{(2k+1)(2m+1)+k-m}{2(2m+1)}$$ 
and hence
\begin{eqnarray*}
  \dfrac{1}{(2k+1)(2m+1)+k-m}\le \dfrac{a}{2k+1}\le \dfrac{1}{(2k+1)(2m+1)+m-k},
\end{eqnarray*}
 \textit{i.e.},
 \begin{eqnarray}\label{rangeofa}
\dfrac{1}{4km+3k+m+1}\le \dfrac{a}{2k+1}\le \dfrac{1}{4km+k+3m+1}.
\end{eqnarray}
Let us define 
\begin{eqnarray}
X_{sn}=\dfrac{s}{2(2m+1)}-\dfrac{n}{2k+1},~\text{for}~s=0,1,\dots,4m+1,~n=0,1,\dots,2k.\label{Xsn}
\end{eqnarray}
Then we have
\begin{align*}
 \dfrac{1}{2(2m+1)}-\dfrac{2k}{2k+1}\le& X_{sn}\le \dfrac{4m+1}{2(2m+1)}-\dfrac{1}{2k+1},~s=\overline{1,4m+1},~n=\overline{1,2k},
 \end{align*}
 which implies that
 \begin{align*}
\dfrac{-8km-2k+1}{8km + 4k + 4m + 2}\le &X_{sn}\le\dfrac{8km+2k-1}{8km + 4k + 4m + 2}.
\end{align*}
Therefore, $-1< X_{sn}<1$ and $X_{sn}\ne 0$ for $s=1,2,\dots,4m+1$ and $n=1,2,\dots,2k$.
Let
\begin{eqnarray}
Y=\dfrac{1}{2a(2m+1)}\label{Y}. 
\end{eqnarray}
We obtain from \eqref{rangeofa} that $$0>2a(2m+1)\ge \dfrac{2(2m+1)(2k+1)}{4km+3k+m+1}=\dfrac{8km+4k+4m+2}{4km+3k+m+1}>1,$$ 
which implies that $0<Y<1$ and $\floor{Y}=-\ceil{-Y}=0$. Here and throughout the paper, 
$\floor{x}$ and $\ceil{x}$ denote the floor and ceiling functions, respectively.

The following lemma provides information about the entries of the symbol matrix $\Phi(0,0)\equiv\Phi_g(0,0)$ for $g=Q_2$.

\begin{lem}\label{lem1}
For $s=1,2,\dots,4m+1$ and $n=1,2,\dots,2k$, we have
$$\Phi_{sn}(0,0)=\begin{cases}
\sum\limits_{l=-1}^0 Q_2\Big(2a(2m+1)(l+X_{sn})\Big), & \text{if}~ 0<X_{sn}<1,\\
\sum\limits_{l=0}^1 Q_2\Big(2a(2m+1)(l+X_{sn})\Big), & \text{if} -1<X_{sn}<0,
\end{cases}$$
where $X_{sn}$ is defined in \eqref{Xsn}.
\end{lem}
\begin{pf}
Recall from \eqref{symbolentries} that
\begin{align*}
  \Phi_{sn}(0,0)=~&\sum\limits_{l\in\mathbb{Z}}Q_2\left(2a(2m+1)l+as-\dfrac{2an(2m+1)}{2k+1}\right)\\
  =~&\sum\limits_{l\in\mathbb{Z}}Q_2\Big(2a(2m+1)(l+X_{sn})\Big).  
\end{align*}
Let $\Lambda=\big\{l\in\mathbb{Z}:~|2a(2m+1)(l+X_{sn})|<1\big\}$. Since $Q_2(x)=(1-|x|)\chi_{[-1,1]}(x)$ is positive in $(-1,1)$ and its support is $[-1,1]$, we have
$$\Phi_{sn}(0,0)=\sum\limits_{l\in\Lambda} Q_2\Big(2a(2m+1)(l+X_{sn})\Big).$$ If $-1<2a(2m+1)l+2a(2m+1)X_{sn}<1$, then
\begin{eqnarray}\label{lem1ineq}
-Y-X_{sn}&<&l<Y-X_{sn},\nonumber\\
\implies -\floor{X_{sn}+Y}&\le& l\le -\ceil{X_{sn}-Y}.
\end{eqnarray}
Since $X_{sn}\notin\mathbb{Z}$, we have $\ceil{X_{sn}}=1+\floor{X_{sn}}$. We know that $\floor{\alpha+\beta}=\floor{\alpha}+\floor{\beta}$ or $\floor{\alpha}+\floor{\beta}+1$ and $\ceil{\alpha+\beta}=\ceil{\alpha}+\ceil{\beta}-1$ or $\ceil{\alpha}+\ceil{\beta}$. Therefore, the following four cases are possible from \eqref{lem1ineq}.\\
\noindent
\underline{\textbf{Case: 1}} When $\floor{X_{sn}+Y}=\floor{X_{sn}}+\floor{Y}$
and $\ceil{X_{sn}-Y}=\ceil{X_{sn}}+\ceil{-Y}$. Since $\ceil{X_{sn}}=1+\floor{X_{sn}}$ and $\floor{Y}=\ceil{-Y}=0$, we have 
$$-\floor{X_{sn}}\le l\le -1-\floor{X_{sn}},$$
which implies that there does not exist any $l$.
\\
\noindent
\underline{\textbf{Case: 2}} When $\floor{X_{sn}+Y}=\floor{X_{sn}}+\floor{Y}+1$ and $\ceil{X_{sn}-Y}=\ceil{X_{sn}}+\ceil{-Y}$, we have 
$$-1-\floor{X_{sn}}\le l\le -1-\floor{X_{sn}}.$$
In this case, $l=-1-\floor{X_{sn}}$.
\\
\noindent
\underline{\textbf{Case: 3}} When $\floor{X_{sn}+Y}=\floor{X_{sn}}+\floor{Y}$ and $\ceil{X_{sn}-Y}=\ceil{X_{sn}}+\ceil{-Y}-1$, we have 
$$-\floor{X_{sn}}\le l\le -\floor{X_{sn}}.$$
In this case, $l=-\floor{X_{sn}}$.
\\
\noindent
\underline{\textbf{Case: 4}} When $\floor{X_{sn}+Y}=\floor{X_{sn}}+\floor{Y}+1$ and $\ceil{X_{sn}-Y}=\ceil{X_{sn}}+\ceil{-Y}-1$, we have 
$$-1-\floor{X_{sn}}\le l\le -\floor{X_{sn}}.$$
In this case, $l=-1-\floor{X_{sn}},-\floor{X_{sn}}$.

From the above four cases we observe that either $\Lambda=\emptyset $ or $\Lambda\subseteq\big\{-1-\floor{X_{sn}},-\floor{X_{sn}}\big\}$. 
When $0<X_{sn}<1$, we have $\floor{X_{sn}}=0$. Therefore, the possible values of $l$ are $-1,0$. When $-1<X_{sn}<0$, we have $\floor{X_{sn}}=-1$ and hence the possible values of $l$ are $l=0,1$.
\end{pf}

Let us define
\begin{eqnarray}
A_{sn}=\Phi_{sn}(0,0)-\Phi_{s,2k+1-n}(0,0),~\text{for}~s = 0,\dots, 4m + 1, ~n = 1, \dots, k.
\end{eqnarray}

\begin{lem}\label{02m+1th0}
    $A_{0n}=A_{2m+1,n}=0$, for $n=1,2,\dots,k$.
\end{lem}
\begin{pf}
It follows from \eqref{rangeofa} that
$$0>2a(2m+1)X_{0n}=\dfrac{-2an(2m+1)}{2k+1}\ge \dfrac{-2k(2m+1)}{4km+k+3m+1}>-1$$
and $$2a(2m+1)(1+X_{0n})=\dfrac{2a(2m+1)(2k+1-n)}{2k+1}\ge \dfrac{4km+2k+4m+2}{4km+k+3m+1}>1.$$
Hence
$$\Phi_{0n}(0,0)=1+2a(2m+1)X_{0n}=1-\dfrac{2an(2m+1)}{2k+1},$$
from Lemma \ref{lem1}.
Again \eqref{rangeofa} yields
$$2a(2m+1)X_{0,2k+1-n}=-\dfrac{2a(2k+1-n)(2m+1)}{2k+1}\le -\dfrac{4km+2k+4m+2}{4km+3k+m+1} <-1,$$
and 
$$0<2a(2m+1)(1+X_{0,2k+1-n})=\dfrac{2an(2m+1)}{2k+1}\le \dfrac{2k(2m+1)}{4km+k+3m+1}<1.$$
Hence 
$$\Phi_{0,2k+1-n}(0,0)=1-2a(2m+1)(1+X_{0,2k+1-n})=1-\dfrac{2an(2m+1)}{2k+1}.$$
Therefore, $A_{0n}=0.$\\
\noindent
   Since 
   $$X_{2m+1,n}=\dfrac{(2k+1-2n)}{2(2k+1)}\ge\dfrac{1}{2(2k+1)}>0,$$
   we obtain from \eqref{rangeofa} that
$$0<2a(2m+1)X_{2m+1,n}=\dfrac{a(2m+1)(2k+1-2n)}{2k+1}\le \dfrac{4km+2k-2m-1}{4km+k+3m+1}<1,$$
because $4km+k+3m+1-(4mk+2k-2m-1)=-k+5m+2>0$.\\
\noindent
Since  $2a(2m+1)(-1+X_{2m+1,n})<0$, using \eqref{rangeofa} we have
\begin{align*}
    2a(2m+1)(-1+X_{2m+1,n})=&~-\dfrac{a(2k+1+2n)(2m+1)}{2k+1}\\
    \le&~-\dfrac{4km+2k+6m+3}{4km+3k+m+1}<-1
\end{align*}
because $4km+3k+m+1-(4km+2k+6m+3)=k-5m-2<0$.
Therefore,
$$\Phi_{2m+1,n}=1-2a(2m+1)X_{2m+1,n}=1-\dfrac{a(2m+1)(2k+1-2n)}{2k+1}.$$
We now compute $\Phi_{2m+1,2k+1-n}$.
Since $X_{2m+1,2k+1-n}=-\dfrac{2k+1-2n}{2(2k+1)}<0$, again \eqref{rangeofa} yields
 \begin{align*}
     0>2a(2m+1)X_{2m+1,2k+1-n}=&\dfrac{a}{2k+1}\big[(2m+1)(2n-2k-1)\big]\\
     \ge &\dfrac{(2m+1)(1-2k)}{4km+k+3m+1}
     =-\dfrac{4mk+2k-2m-1}{4km+k+3m+1}>-1
 \end{align*}
and
$$2a(2m+1)(1+X_{2m+1,2k+1-n})=\dfrac{a(2k+1+n)(2m+1)}{2k+1}\ge \dfrac{(2k+2)(2m+1)}{4mk+k+3m+1}>1.$$ 
Therefore,
$$\Phi_{2m+1,2k+1-n}=1+2a(2m+1)X_{2m+1,2k+1-n}=1-\dfrac{a(2m+1)(2k+1-2n)}{2k+1}.$$
Hence $A_{2m+1,n}=0.$
\end{pf}

\begin{lem}\label{reflection}
    For $s=1,2,\dots,2m$ and  $n=1,2,\dots,k$, we have
   $$A_{s,n}=-A_{4m+2-s,n}.$$
\end{lem}
\begin{pf}
Using the symmetric property of $Q_2$, we have
    \begin{align}
        \Phi_{4m+2-s,n}(0,0)=~&\sum\limits_{l=-1}^{1}Q_2\bigg(2a(2m + 1)(l + X_{4m+2-s,n})\bigg)\nonumber\\ =~&\sum\limits_{l=-1}^{1}Q_2\bigg(2a(2m + 1)l-as+\dfrac{2a(2k+1-n)(2m+1)}{2k+1}\bigg)\nonumber\\
=~&\sum\limits_{l=-1}^{1}Q_2\bigg(2a(2m + 1)l+as-\dfrac{2a(2k+1-n)(2m+1)}{2k+1}\bigg)\nonumber\\
=~&\Phi_{s,2k+1-n}(0,0).\nonumber
    \end{align}
    Similarly, we can show that $\Phi_{4m+2-s,2k+1-n}(0,0)=\Phi_{sn}(0,0).$
    Therefore,
    \begin{align*}
        A_{4m+2-s,n}=~&\Phi_{4m+2-s,n}-\Phi_{4m+2-s,2k+1-n}(0,0)\\
        =~&\Phi_{s,2k+1-n}(0,0)-\Phi_{sn}(0,0)=- A_{sn}.
     \end{align*}
\end{pf}

In the following sequence of Lemmas, we assume that $k$ is even.
\begin{lem}\label{lem2}
    For $s = 1, 2,\dots, m$ and $n = 1, 2, \dots, k/2$, we have

\begin{align}\label{Asnlem2}
A_{sn}=\begin{cases}
{2n}/{b}, & \text{if}~ 0<X_{sn}<Y,\\
2as, & \text{if} -Y<X_{sn}<0.
\end{cases} 
    \end{align}
\end{lem}
\begin{pf}
When $0<X_{sn}<1$, it follows from \eqref{rangeofa} that
\begin{align}
0<~&2a(2m+1)X_{sn}
=\dfrac{a}{2k+1}\big[s(2k+1)-2n(2m+1)\big]\nonumber\\
\le~& \dfrac{s(2k+1)-2n(2m+1)}{4km+k+3m+1}
     \le\dfrac{m(2k+1)-2(2m+1)}{4km+k+3m+1}\nonumber\\
     =~&\dfrac{2km-3m-2}{4km+k+3m+1}<1.\nonumber
    \end{align}
Since $2a(2m+1)(-1+X_{sn})<0$,  we obtain from \eqref{rangeofa} that
 \begin{align} 
2a(2m+1)(-1+X_{sn})=~&\dfrac{a}{2k+1}\big[-2(2m+1)(2k+1)+s(2k+1)-2n(2m+1)\big]\nonumber \\
 \le~&  \dfrac{-2(2m+1)(2k+1)+s(2k+1)-2n(2m+1)}{4km+3k+m+1}\nonumber\\
 \le~& \dfrac{-2(2m+1)(2k+1)+m(2k+1)-2(2m+1)}{4km+3k+m+1}\nonumber\\
 =~&-\dfrac{6km+4k+7m+4}{4km+3k+m+1}<-1.\nonumber
    \end{align}
    Consequently, we obtain from Lemma \ref{lem1} that
\begin{align}\label{phisnlem2case1}
        \Phi_{sn}(0,0)
        =1-2a(2m+1)X_{sn}~\text{if}~0<2a(2m+1)X_{sn}<1.
    \end{align}
 When $-1<X_{sn}<0$, again using \eqref{rangeofa} we have 
    \begin{align}
        0>~&2a(2m+1)X_{sn}
        =\dfrac{a}{2k+1}\big[s(2k+1)-2n(2m+1)\big]\nonumber\\
        \ge~&\dfrac{(2k+1)-k(2m+1)}{4km+k+3m+1}
        =\dfrac{-2km+k+1}{4km+k+3m+1}>-1\nonumber.
    \end{align}
Since 
$2a(2m+1)(1+X_{sn})>0$, using \eqref{rangeofa} we get
\begin{align}
2a(2m+1)(1+X_{sn})\nonumber=~&\dfrac{a}{2k+1}\Big[2(2m+1)(2k+1)+s(2k+1)-2n(2m+1)\Big]\nonumber\\
\ge~&\dfrac{2(2m+1)(2k+1)+(2k+1)-k(2m+1)}{4km+3k+m+1}\nonumber\\
=~&\dfrac{6km+5k+4m+3}{4km+3k+m+1}>1\nonumber.
\end{align}
Consequently, we obtain from Lemma \ref{lem1} that    \begin{align}\label{phisnlem2case2}
        \Phi_{sn}(0,0)
=1+2a(2m+1)X_{sn}~\text{if}~-1<2a(2m+1)X_{sn}<0.
\end{align}
We now compute $\Phi_{s,2k+1-n}(0,0).$
Since
\begin{align}
    X_{s,2k+1-n}=~&\dfrac{s(2k+1)-2(2k+1-n)(2m+1)}{2(2m+1)(2k+1)}\nonumber\\
    \le~&\dfrac{m(2k+1)-2(2k+1)(2m+1)+k(2m+1)}{2(2m+1)(2k+1)}\nonumber\\
    =~&-\dfrac{4km+3k+3m+2}{8km+4k+4m+2}<0\nonumber,
\end{align}
we obtain from \eqref{rangeofa} that
\begin{align}
2a(2m+1)X_{s,2k+1-n}\le\dfrac{-(4km+3k+3m+2)}{4km+3k+m+1}<-1\nonumber.
\end{align}
Again \eqref{rangeofa} yields  
\begin{align}
0<~&2a(2m+1)(1+X_{s,2k+1-n})=\dfrac{a}{2k+1}\big[s(2k+1)+2n(2m+1)\big]\nonumber\\
    \le~&\dfrac{m(2k + 1)+k(2m + 1)}{4km +k+ 3m + 1}    =\dfrac{4km+k+m}{4km +k+ 3m + 1}<1\nonumber.
\end{align}
Therefore, we obtain from Lemma \ref{lem1} that
\begin{align}\label{phisnl=1lem2}
    \Phi_{s,2k+1-n}(0,0)=1-2a(2m+1)-2a(2m+1)X_{s,2k+1-n}.
\end{align}
The desired result  now follows from \eqref{phisnlem2case1}, \eqref{phisnlem2case2}, and \eqref{phisnl=1lem2}.
\end{pf}

\begin{lem}\label{lem3}
    For $s = 1, 2,\dots, m$ and $n = {k}/{2}+1, \dots, k$, we have 
\begin{align}\label{Asnlem3}
A_{sn}=
\begin{cases}
{(2k+1-2n)}/{b}, & \text{if}~ -Y<X_{s,2k+1-n}<0,\\
2as, & \text{if} -1<X_{s,2k+1-n}<Y-1,\\
-1+as+{(2k+1-n)}/{b}, &\text{if}~-Y<X_{s,2k+1-n}<Y-1,\\
1+as-{n}/{b},&\text{if }Y-1\le X_{s,2k+1-n}\le -Y.
\end{cases}    
\end{align}
\end{lem}
\begin{pf}
Since 
\begin{align*}
    X_{sn}=\dfrac{s(2k+1)-2n(2m+1)}{2(2k+1)(2m+1)}    \le~&\dfrac{m(2k+1)-(k+2)(2m+1)}{2(2k+1)(2m+1)}\\
    =~&\dfrac{-3m-k-2}{2(2k+1)(2m+1)}<0,
    \end{align*}
\eqref{rangeofa} yields 
\begin{align}
0>2a(2m+1)X_{sn}\ge~&\dfrac{s(2k + 1)- 2n(2m + 1)}{4km+k+3m+1}\nonumber\\
\ge~&\dfrac{(2k + 1)- 2k(2m + 1)}{4km+k+3m+1},\nonumber\\
=~&\dfrac{-4km+1}{4km+k+3m+1}>-1.
\end{align}
Since $2a(2m+1)(1+X_{sn})>0$, it follows from \eqref{rangeofa} that
\begin{align}
2a(2m+1)(1+X_{sn})=~&\dfrac{a}{2k+1}\Big[2(2m+1)(2k+1)+s(2k+1)-2n(2m+1)\Big]\nonumber\\
\ge~&\dfrac{2(2k+1)(2m+1)+(2k + 1) - 2k(2m + 1)}{4km + 3k + m + 1}\nonumber\\
=~&\dfrac{4km+4k+4m+3}{4km + 3k + m + 1}>1.\nonumber
\end{align}
Hence 
\begin{align}\label{phisnlem3first}
    \Phi_{sn}(0,0)=1+2a(2m+1)X_{sn}.
\end{align}
We now compute $\Phi_{s,2k+1-n}(0,0)$. Since
\begin{align*}
    X_{s,2k+1-n}=~&\dfrac{s(2k+1)-2(2k+1-n)(2m+1)}{2(2m+1)(2k+1)}\\
    \le~&\dfrac{m(2k+1)-2(k+1)(2m+1)}{2(2m+1)(2k+1)}
    =\dfrac{-2km-2k-3m-2}{2(2m+1)(2k+1)}<0,
\end{align*} 
it follows from \eqref{rangeofa} that 
\begin{align}\label{2a(2m+1)ineq}
   \dfrac{-6km-k+1}{4km+k+3m+1} \le~ 2a(2m+1)X_{s,2k+1-n}\le~\dfrac{-2km-2k-3m-2}{4km+3k+m+1}<0.
\end{align}
Since $$ \dfrac{-6km-k+1}{4km+k+3m+1} <-1 \text{ for }m>1 \text{ and } -1<\dfrac{-2km-2k-3m-2}{4km+3k+m+1},$$
it follows from \eqref{2a(2m+1)ineq} that either $2a(2m+1)X_{s,2k+1-n}\le-1$ or $-1<2a(2m+1)X_{s,2k+1-n}<0$.

Similarly, we can show that $2a(2m+1)(1+X_{s,2k+1-n})>0$
and hence
\begin{align}\label{2a(2m+1)l=1}
\dfrac{2km + 3k + 4m + 3}{4km + 3k + m + 1}\le~2a(2m+1)(1+X_{s,2k+1-n})\le~\dfrac{6km+2k+m}{4km+k+3m+1}
\end{align}
Since 
$$\dfrac{2km + 3k + 4m + 3}{4km + 3k + m + 1}< 1 \text{ for $m>1$ and } 1<\dfrac{6km+2k+m}{4km+k+3m+1},$$
it follows from \eqref{2a(2m+1)l=1} that either $0<2a(2m+1)(1+X_{s,2k+1-n})<1$ or $1\le2a(2m+1)(1+X_{s,2k+1-n})$. Therefore, the following four cases occur.\\
\noindent
\underline{\textbf{Case:1}} When $-1<2a(2m+1)X_{s,2k+1-n}<0$ and $1\le2a(2m+1)(1+X_{s,2k+1-n})$, we have 
\begin{align}\label{phisnlem3secondcase1}
\Phi_{s,2k+1-n}(0,0)=
1+2a(2m+1)X_{s,2k+1-n}.
\end{align}\\
\noindent
\underline{\textbf{Case:2}}
When $2a(2m+1)X_{s,2k+1-n}\le-1$ and $0<2a(2m+1)(1+X_{s,2k+1-n})<1$, we have
\begin{align}\label{phisnlem3secondcase2}
\Phi_{s,2k+1-n}(0,0)&=1-2a(2m+1)-2a(2m+1)X_{s,2k+1-n}\nonumber\\
&=1-as-\dfrac{2na(2m+1)}{2k+1}.
\end{align}
\\
\noindent
\underline{\textbf{Case:3}}
When $-1<2a(2m+1)X_{s,2k+1-n}<0$ and $0<2a(2m+1)(1+X_{s,2k+1-n})<1$, we have 
\begin{align}\label{phisnlem3secondcase3}
\Phi_{s,2k+1-n}(0,0)
=&\Big(1+2a(2m+1)X_{s,2k+1-n}\Big)+\Big(1-2a(2m+1)(1+X_{s,2k+1-n})\Big)\nonumber\\
=&2-2a(2m+1).
\end{align}
\\
\noindent
\underline{\textbf{Case:4}}
When $2a(2m+1)X_{s,2k+1-n}\le-1$ and $1\le2a(2m+1)(1+X_{s,2k+1-n})$, we have
\begin{align}\label{phisnlem3secondcase4}
\Phi_{s,2k+1-n}(0,0)=0.
\end{align}
 Hence we can conclude our result using
 \eqref{phisnlem3first}, \eqref{phisnlem3secondcase1}, \eqref{phisnlem3secondcase3}, \eqref{phisnlem3secondcase2}, and \eqref{phisnlem3secondcase4}.
\end{pf}
\begin{lem}\label{lem4}
    For $s = m+1, m+2,\dots, 2m$ and $n = 1,2, \dots, k/2$, we have 
\begin{align}\label{Asnlem4}
A_{sn}=\begin{cases}
2a(2m+1-s), & \text{if}~ -Y<X_{s,2k+1-n}<0,\\
{2n}/{b}, & \text{if} -1<X_{s,2k+1-n}<Y-1,\\
-1-as+{(2k+1+n)}/{b}, &\text{if}~-Y<X_{s,2k+1-n}<Y-1,\\
1-as+{n}/{b},&\text{if }Y-1\le X_{s,2k+1-n}\le -Y.
\end{cases}    
\end{align}
\end{lem}
\begin{pf}
    Since 
    \begin{align*}
    X_{sn}=\dfrac{s(2k+1)-2n(2m+1)}{2(2k+1)(2m+1)}
    \ge~&\dfrac{(m+1)(2k+1)-k(2m+1)}{2(2k+1)(2m+1)}\\
    =~&\dfrac{m+k+1}{2(2k+1)(2m+1)}>0,
    \end{align*}
    using \eqref{rangeofa} we have 
\begin{align}\label{lem4l=0first}
0<~&2a(2m+1)X_{sn}\le\dfrac{s(2k + 1)- 2n(2m + 1)}{4km+k+3m+1}\nonumber\\
\le~&\dfrac{2m(2k + 1)- 2(2m + 1)}{4km+k+3m+1}
=\dfrac{4km-2m-2}{4km+k+3m+1}<1.
\end{align}
Since $2a(2m+1)(-1+X_{sn})<0,$
\eqref{rangeofa} yields
\begin{align*}
2a(2m+1)(-1+X_{sn})
=~&\dfrac{a}{2k+1}\big[-2(2m+1)(2k+1)+s(2k+1)-2n(2m+1)\big]\nonumber\\
\le~&\dfrac{a}{2k+1}\big[-2(2m+1)(2k+1)+2m(2k+1)-2(2m+1)\big]\nonumber\\
\le~&-\dfrac{4km+4k+6m+4}{4km+3k+m+1}<-1.\nonumber
\end{align*}
Consequently, we obtain from Lemma \ref{lem1} that
\begin{align}\label{phisnlem4first}
  \Phi_{sn}(0,0)=1-2a(2m+1)X_{sn}.
\end{align}
Now we compute $\Phi_{s,2k+1-n}(0,0)$. Since
\begin{align*}
    X_{s,2k+1-n}=~&\dfrac{s(2k+1)-2(2k+1-n)(2m+1)}{2(2m+1)(2k+1)}\\
    \le~&\dfrac{2m(2k+1)-(3k+2)(2m+1)}{2(2m+1)(2k+1)}=\dfrac{-2km-3k-2m-2}{2(2m+1)(2k+1)}<0,
\end{align*} 
using \eqref{rangeofa} we have 
\begin{align}\label{2a(2m+1)lem4l=0}
   \dfrac{-6km-2k+m+1}{4km+k+3m+1} \le~ 2a(2m+1)X_{s,2k+1-n}\le~\dfrac{-2km-3k-2m-2}{4km+3k+m+1}<0.
\end{align}
Since
\begin{align*}
    \dfrac{-6km-2k+m+1}{4km+k+3m+1}<-1 \text{ for } m>1 \text{ and }-1<\dfrac{-2km-3k-2m-2}{4km+3k+m+1},
\end{align*}
it follows from \eqref{2a(2m+1)lem4l=0} that either $2a(2m+1)X_{s,2k+1-n}\le-1$ or $-1<2a(2m+1)X_{s,2k+1-n}<0$.

Similarly, we can show that $2a(2m+1)(1+X_{s,2k+1-n})>0$
and hence 
\begin{align}\label{2a(2m+1)lem4l=1}
\dfrac{2km + 2k + 3m + 2}{4km + 3k + m + 1}\le~2a(2m+1)(1+X_{s,2k+1-n})\le~\dfrac{6km+k+2m}{4km+k+3m+1}.
\end{align}
Since
\begin{align*}
    \dfrac{2km + 2k + 5m + 3}{4km + 3k + m + 1}<1\text{ for }m>1\text{ and } 1<\dfrac{6km+k+2m}{4km+k+3m+1},
\end{align*}
it follows from \eqref{2a(2m+1)lem4l=1} that either $0<2a(2m+1)(1+X_{s,2k+1-n})<1$ or $1\le2a(2m+1)(1+X_{s,2k+1-n})$.
Therefore, the following four cases occur.\\
\noindent
\underline{\textbf{Case: 1}}
When $-1<2a(2m+1)X_{s,2k+1-n}<0$ and $1\le2a(2m+1)(1+X_{s,2k+1-n})$, we have 
\begin{align}\label{phisnlem4secondcase1}
\Phi_{s,2k+1-n}(0,0)=1+2a(2m+1)X_{s,2k+1-n}.
\end{align}
\\
\noindent
\underline{\textbf{Case: 2}}
When $-1<2a(2m+1)X_{s,2k+1-n}<0$ and $0<2a(2m+1)(1+X_{s,2k+1-n})<1$, we have 
\begin{align}\label{phisnlem4secondcase3}
\Phi_{s,2k+1-n}(0,0)&=1+2a(2m+1)X_{s,2k+1-n}+1-2a(2m+1)(1+X_{s,2k+1-n})\nonumber\\
&=2-2a(2m+1).
\end{align}
\\
\noindent
\underline{\textbf{Case: 3}}
When $2a(2m+1)X_{s,2k+1-n}\le-1$ and $0<2a(2m+1)(1+X_{s,2k+1-n})<1$, we have
\begin{align}\label{phisnlem4secondcase2}
\Phi_{s,2k+1-n}(0,0)&=1-2a(2m+1)(1+X_{s,2k+1-n}).
\end{align}
\\
\noindent
\underline{\textbf{Case: 4}}
When $2a(2m+1)X_{s,2k+1-n}\le-1$ and $1\le2a(2m+1)(1+X_{s,2k+1-n})$, we have
\begin{align}\label{phisnlem4secondcase4}
\Phi_{s,2k+1-n}(0,0)=0.
\end{align}
Hence we can conclude our result using \eqref{phisnlem4first}, \eqref{phisnlem4secondcase1}, \eqref{phisnlem4secondcase3}, \eqref{phisnlem4secondcase2}, and \eqref{phisnlem4secondcase4}.
\end{pf}
\begin{rem}\label{remlem4}
   By replacing $s$ by $2m+1-s$, \eqref{Asnlem4} can be written as $A_{2m+1-s,n}$
\begin{align}\label{Asnlem4new}
=\begin{cases}
2as, & \text{if}~ -Y<X_{2m+1-s,2k+1-n}<0,\\
{2n}/{b}, & \text{if} -1<X_{2m+1-s,2k+1-n}<Y-1,\\
-1-a(2m+1-s)+{(2k+1+n)}/{b}, &\text{if}~-Y<X_{2m+1-s,2k+1-n}<Y-1,\\
1-a(2m+1-s)+{n}/{b}, &\text{if }Y -1 \le X_{2m+1-s,2k+1-n} \le -Y,
\end{cases}    
\end{align}
for $s = 1, 2,\dots, m$ and $n = 1, \dots, k/2$.
\end{rem}

\begin{lem}\label{lem5}
    For $s = m+1, m+2,\dots, 2m$ and $n = k/2+1, \dots, k$, we have
\begin{align}\label{Asnlem5}
    A_{sn}=\begin{cases}
2a(2m+1-s), & \text{if}~ 0<X_{sn}<Y,\\
{(2k+1-2n)}/{b}, & \text{if} -Y<X_{sn}<0.
\end{cases}
\end{align}
\end{lem}
\begin{pf}
    When $0<X_{sn}<1$, using \eqref{rangeofa} we have 
    \begin{align}
     0<~&2a(2m+1)X_{sn}
     =\dfrac{a}{2k+1}\big[s(2k+1)-2n(2m+1)\big]\nonumber\\
     \le~& \dfrac{s(2k+1)-2n(2m+1)}{4km+k+3m+1}
     \le\dfrac{2m(2k+1)-(k+2)(2m+1)}{4km+k+3m+1}\nonumber\\
     =~&\dfrac{2km-k-2m-2}{4km+k+3m+1}<1.\nonumber
    \end{align}
    Since $2a(2m+1)(-1+X_{sn})<0$, using \eqref{rangeofa} we get $2a(2m+1)(-1+X_{sn})$
    \begin{align*}
=~& \dfrac{a}{2k+1}\big[-2(2m+1)(2k+1)+s(2k+1)-2n(2m+1)\big]\\
\le~& \dfrac{a}{2k+1}\big[-2(2m+1)(2k+1)+2m(2k+1)-(k+2)(2m+1)\big]\\
\le~&-\dfrac{6km + 5k + 6m + 4}{4km+3k+m+1}<-1.
    \end{align*}
 Therefore, 
    \begin{align}\label{phisnlem5case1}
        \Phi_{sn}(0,0)=1-2a(2m+1)X_{sn},~\text{if}~0<2a(2m+1)X_{sn}<1.
    \end{align}
    When $-1<X_{sn}<0$, using \eqref{rangeofa} we have 
    \begin{align}
       0> 2a(2m+1)X_{sn}
        =~&\dfrac{a}{2k+1}\big[s(2k+1)-2n(2m+1)\big]\nonumber\\
        \ge~&\dfrac{(m+1)(2k+1)-2k(2m+1)}{4km+k+3m+1}\nonumber\\
        =~&\dfrac{-2km+m+1}{4km+k+3m+1}>-1,\nonumber
    \end{align}
    Since $2a(2m+1)(1+X_{sn})>0$, using \eqref{rangeofa} we have
    \begin{align}
2a(2m+1)(1+X_{sn})
        =~&\dfrac{a}{2k+1}\big[2(2m+1)(2k+1)+s(2k+1)-2n(2m+1)\big]\nonumber\\
        \ge~&\dfrac{2(2m+1)(2k+1)+(m+1)(2k+1)-2k(2m+1)}{4km+3k+m+1}\nonumber\\
\ge~&\dfrac{6km+4k+5m+3}{4km+3k+m+1}>1.\nonumber
    \end{align}
Therefore, 
    \begin{align}\label{phisnlem5case2}
        \Phi_{sn}(0,0)=1+2a(2m+1)X_{sn},~\text{if}~-1<2a(2m+1)X_{sn}<0.
    \end{align}
    We now compute $\Phi_{s,2k+1-n}(0,0).$
Since
\begin{align}
    X_{s,2k+1-n}=~&\dfrac{s(2k+1)-2(2k+1-n)(2m+1)}{2(2m+1)(2k+1)}\nonumber\\
    \le~&\dfrac{2m(2k+1)-2(2k+1)(2m+1)+2k(2m+1)}{2(2m+1)(2k+1)}\nonumber\\
    =~&\dfrac{-(2k+2m+2)}{2(2m+1)(2k+1)}<0\nonumber,
\end{align}
using \eqref{rangeofa} we have 
\begin{align*}
0>2a(2m+1)X_{s,2k+1-n}=&\dfrac{a}{2k+1}\big[s(2k+1)-2(2k+1-n)(2m+1)\big]\\
    \ge&\dfrac{(m+1)(2k+1)-2(2k+1)(2m+1)+(k+2)(2m+1)}{4km+k+3m+1}\\
    \ge&-\dfrac{4km+k-m-1}{4km+k+3m+1}>-1.
\end{align*}
Since $2a(2m+1)(1+X_{s,2k+1-n})>0$, using \eqref{rangeofa} we have
\begin{align}
2a(2m+1)(1+X_{s,2k+1-n})
=~&\dfrac{a}{2k+1}\big[s(2k+1)+2n(2m+1)\big]\nonumber\\
    \ge~&\dfrac{(m+1)(2k + 1) + (k+2)(2m + 1)}{4km +k+ 3m + 1}\nonumber\\
    =~&\dfrac{4km+3k+5m+3}{4km + k+3m + 1}>1\nonumber.
\end{align}
Therefore,
\begin{align}\label{phisnl=1lem5}
    \Phi_{s,2k+1-n}(0,0)= 1-2a(2m+1)X_{s,2k+1-n}
\end{align}
from Lemma \ref{lem1}.
From \eqref{phisnlem5case1}, \eqref{phisnlem5case2}, and \eqref{phisnl=1lem5} we obtain our result.
\end{pf}
\begin{rem}\label{remlem5}
     By replacing $s$ by $2m+1-s$, \eqref{Asnlem5} can be written as 
\begin{align}\label{Asnlem5new}
    A_{2m+1-s,n}=\begin{cases}
2as, & \text{if}~ 0<X_{2m+1-s,n}<Y,\\
{(2k+1-2n)}/{b}, & \text{if} -Y<X_{2m+1-s,n}<0.
\end{cases}
\end{align}
for $s = 1, 2,\dots, m$ and $n = k/2+1, \dots, k$.
\end{rem}

For each $\mu\in\mathbb{N}$, let us define $J_\mu=\{1,2,\dots,\mu\}$. Let 
\begin{align}
S_{1a}=&\bigg\{s\in J_m: A_{2m+1-s,n}=1-a(2m+1-s)+\dfrac{n}{b},\text{ for some }n\in J_{\tfrac{k}{2}}\bigg\},\\
S_{2a}=&\bigg\{s\in J_m: A_{s,n}=1+as-\dfrac{n}{b},\text{ for some }n\in \tfrac{k}{2}+J_{\tfrac{k}{2}}\bigg\},\\
S_{3a}=&\bigg\{s\in J_m: A_{2m+1-s,n}=-1-a(2m+1-s)+\dfrac{2k+1+n}{b},\text{ for some }n\in J_{\tfrac{k}{2}}\bigg\},
\end{align}
and
\begin{align}
S_{4a}=&\bigg\{s\in J_m: A_{s,n}=-1+as+\dfrac{2k+1-n}{b},\text{ for some }n\in \tfrac{k}{2}+J_{\tfrac{k}{2}}\bigg\}.
    \end{align}    
    
Now we observe from Remark 
    \ref{remlem4} (or Case: $4$ in Lemma $\ref{lem4}$) that 
    $$A_{2m+1-s,n}=1-a(2m+1-s)+{n}/{b}$$ 
    if and only if 
    $$2a(2m + 1)X_{2m+1-s,2k+1-n} \le -1\text{ and }1\le 2a(2m + 1)(1+X_{2m+1-s,2k+1-n}).$$
    The above conditions are equivalent to 
    \begin{align}\label{ineqsoddcase}
        1\le a(2m+1-s)+\dfrac{2a(2m+1)n}{2k+1}\le 2a(2m+1)-1.
    \end{align}
     We write \eqref{ineqsoddcase} as
    \begin{align}\label{rangeofslem4last}
        \ceil{W_n-V}\le~ s~\le \floor{W_n+V},
    \end{align}
    where $W_n=\dfrac{2(2m+1)n}{2k+1}$ and $V=2m+1-\dfrac{1}{a}$. We obtain from \eqref{rangeofa} that
    \begin{align*}
        2m+1-\dfrac{4km+3k+m+1}{2k+1}~\le~V~\le~ 
        2m+1-\dfrac{4km+k+3m+1}{2k+1},  
    \end{align*}
    which implies that
    \begin{align*}
        -1<-\dfrac{k-m}{2k+1}~\le~& V~\le~ \dfrac{k-m}{2k+1}<1.
    \end{align*}
    
When $0\le V<1$, we have $\floor{V}=-\ceil{-V}=0$. Arguing as in Lemma $\ref{lem1}$, we can conclude that for each $n\in J_{{k}/{2}}$, the maximum possible values of $s\in S_{1a}$ are $\floor{W_n}$ and $ 1+\floor{W_n}.$ When $-1<V<0$, we can not find any possible values of $s$ such that \eqref{rangeofslem4last} holds. In this case, $S_{1a}=\emptyset$. When $0\le V<1$, using a similar argument we can show that the maximum possible values of $s\in S_{2a}$ are $2m+1-\floor{W_n}$ and $2m-\floor{W_n}$ for each $n\in k/2+J_{k/2}$ and when $-1<V<0,$ $S_{2a}=\emptyset.$ Now it is clear that $S_{1a}$ and $S_{2a}$ are disjoint.

Similarly, we can possibly find $s\in S_{ja}$, for $j=3,4$ only if $-1<V<0$. On the other hand, if $0\le V<1$, then $S_{3a}=S_{4a}=\emptyset$. It is also easy to check that $S_{3a}$ and $S_{4a}$ are disjoint. Therefore, $S_{1a}\cup S_{2a}$ is mutually disjoint with $S_{3a}\cup S_{4a}$. Hence it is clear that when $0\le V<1$, then $S_{3a}\cup S_{4a}=\emptyset$ and when $-1< V<0$, then $S_{1a}\cup S_{2a}=\emptyset$.\\
\noindent
Now using \eqref{ineqsoddcase} we can rewrite
\begin{align}\label{conditionS1}
   S_{1a}=\bigg\{s\in J_m: \tfrac{1}{a}-W_n\le 2m+1-s\le 2(2m+1)-W_n-\tfrac{1}{a},\text{ for some }n\in J_{\tfrac{k}{2}}\bigg\}.
\end{align}
Similarly, using Case: $4$ in the proof of Lemma $\ref{lem3}$ we can rewrite 
\begin{align}
 S_{2a}=&\bigg\{s\in J_m: \tfrac{1}{a}-W_n\le s\le 2(2m+1)-W_n-\tfrac{1}{a},\text{ for some }n\in\tfrac{k}{2}+J_{\tfrac{k}{2}}\bigg\}\label{conditionoldS2}\\
    =&\bigg\{s\in J_m: \tfrac{1}{a}-W_{n+\tfrac{k}{2}}\le s\le 2(2m+1)-W_{n+\tfrac{k}{2}}-\tfrac{1}{a},\text{ for some }n\in J_{\tfrac{k}{2}}\bigg\}\label{conditionS2}.
\end{align}
Let $\alpha=2m+1-\dfrac{k-m}{2k+1}$. Now we obtain from \eqref{conditionS1} and \eqref{conditionoldS2} that
\begin{align}\label{newS1alpha}
        S_{1,1/\alpha}=&\bigg\{s\in J_m: \tfrac{2n(2m+1)-(k-m)}{2k+1}\le s\le \tfrac{2n(2m+1)+k-m}{2k+1},\text{ for some } n\in J_{\tfrac{k}{2}}\bigg\}
    \end{align}
    and $S_{2,1/\alpha}$
    \begin{align}\label{newS2alpha}
        =\bigg\{s\in J_m: \tfrac{2n(2m+1)-(k-m)}{2k+1}\le 2m+1-s\le \tfrac{2n(2m+1)+k-m}{2k+1},\text{ for some } n\in J_{\tfrac{k}{2}+{\tfrac{k}{2}}}\bigg\}.
    \end{align}
  Let us define
  \begin{align}
    I_{11}=&\bigg\{s\in J_{m}: s=\dfrac{2n(2m+1)+k-m}{2k+1},\text{ for some }n\in J_{\tfrac{k}{2}}\bigg\},\label{L11}\\
    I_{12}=&\bigg\{s\in J_{m}: 2m+1-s=\dfrac{2n(2m+1)+k-m}{2k+1},\text{ for some }n\in \tfrac{k}{2}+J_{\tfrac{k}{2}}\bigg\}\label{L12},\\
    I_{21}=&\bigg\{s\in J_{m}: s=\dfrac{2n(2m+1)-(k-m)}{2k+1},\text{ for some }n\in J_{\tfrac{k}{2}}\bigg\}\label{L21},
\end{align}
and
\begin{align}\label{L22}
     I_{22}=\bigg\{s\in J_{m}: 2m+1-s=\dfrac{2n(2m+1)-(k-m)}{2k+1},\text{ for some }n\in \tfrac{k}{2}+J_{\tfrac{k}{2}}\bigg\}.
\end{align}
\begin{lem}$($\cite{cimnt}, p-$32)$. \label{uniqueness}
The congruence $ax\equiv b\pmod{r}$ has a unique solution in $\mathbb{Z}_{r}$ if and only if $\gcd{(a, r)}=1$.
\end{lem}
\begin{lem}\label{I11unionI2}
$I_{21}=\{m\}$ and $I_{11}=I_{12}= I_{22}=\emptyset$.
\end{lem}
\begin{pf} By choosing $n=k/2$, $m\in I_{21}$. Since $\gcd{(2(2m+1),2k+1)}=1$, we can conclude that $I_{21}=\{m\}$ from Lemma \ref{uniqueness}. Next we prove that $I_{11}=\emptyset$. If possible, $I_{11}$ is nonempty. Let $s_0\in I_{11}$. Then there exists $n_0\in J_{k/2}$ such that $\tfrac{2n_0(2m+1)+k-m}{2k+1}=s_0\in J_m$ which implies that $2n_0(2m+1)\equiv m-k\pmod{2k+1}$. Since $\gcd{(2(2m+1),2k+1)}=1$, $n_0\in J_{k/2}$ must be unique from Lemma \ref{uniqueness}.
Now choose $n_1=2k+1-n_0$. Then
$$\dfrac{2(2k+1-n_0)(2m+1)-(k-m)}{2k+1}=2(2m+1)-s_0.$$
Therefore, $n_1=2k+1-n_0\in \{\tfrac{3k}{2}+1,\dots, 2k\}$ is a solution of $2n_1(2m+1)\equiv k-m\pmod{2k+1}$. But $k/2\ne n_1$ is also a solution of $2n(2m+1)\equiv k-m \pmod{2k+1}$ which contradicts the unique solution in Lemma \ref{uniqueness}. Similarly, we can show that $I_{12}=I_{22}=\emptyset$. 
\end{pf}
\begin{lem}\label{cardslem}
Let $S_a=S_{1a}\cup S_{2a}$ and $T_a=S_{3a}\cup S_{4a}$.
    \noindent
    \begin{enumerate}
        \item [$(i)$] If $2m+1-\dfrac{k-m}{2k+1}\le \dfrac{1}{a}\le ~2m+1$, then $\#S_{a}\le k-m$.
        \item[$(ii)$] If $2m+1< \dfrac{1}{a}\le~ 2m+1+\dfrac{k-m}{2k+1}$, then $\#T_a\le k-m$.
    \end{enumerate}
    Here $\# I$ denotes the cardinality of a finite set $I$.
\end{lem}

\begin{pf}
$(i)$
Since $S_{1a}$ and $S_{2a}$ are disjoint sets and both are monotonically decreasing family of sets with respect to $1/a$, we have 
$\#S_{a}\le \#S_{{1}/{\alpha}}.$ 
We now show that $\#S_{{1}/{\alpha}}=k-m$.
From \eqref{conditionS2}, we have 
    \begin{align}
        S_{2,1/\alpha}=&\bigg\{s\in J_m: W_{n+\tfrac{k}{2}}-\tfrac{k-m}{2k+1}\le 2m+1-s\le W_{n+\tfrac{k}{2}}+\tfrac{k-m}{2k+1},\text{ for some } n\in J_{\tfrac{k}{2}}\bigg\}\nonumber
    \end{align}
Since
\begin{align*}
 W_{n+\tfrac{k}{2}}-\dfrac{k-m}{2k+1}=~&\dfrac{(k+2n)(2m+1)}{2k+1}-\dfrac{k-m}{2k+1}=\dfrac{2km+4mn+2n+m}{2k+1}\\
 =~&m+\dfrac{2n(2m+1)}{(2k+1)}=m+2n-\dfrac{4n(k-m)}{2k+1},
\end{align*}
we have
\begin{align}\label{gauss1}
        S_{2,1/\alpha}=\bigg\{s\in J_m:  \tfrac{(4n-2)(k-m)}{2k+1}\le 2n+s-m-1\le \tfrac{4n(k-m)}{2k+1},\text{ for some } n\in J_{\tfrac{k}{2}}\bigg\}.
\end{align}
Similarly, we can show from \eqref{conditionS1} that
\begin{align}\label{gauss2}
    S_{1,1/\alpha}
        =&\bigg\{s\in J_m:  \tfrac{(4n-1)(k-m)}{2k+1}\le 2n-s\le \tfrac{(4n+1)(k-m)}{2k+1},\text{ for some } n\in J_{\tfrac{k}{2}}\bigg\}.
\end{align}
It is clear that $\#S_{2,1/\alpha}$ is the number of lattice points inside and on the region bounded by the four lines $L_1:y_1(n)= \tfrac{4n(k-m)}{2k+1}$, $L_2:y_2(n)= \tfrac{(4n-2)(k-m)}{2k+1}$, $L_3: n=1$, and $L_4: n=k/2$ (see Figure \ref{cardfig}). Since $\gcd(2(2m+1),2k+1)=\gcd(k-m,2k+1)=1$ and $I_{12}=\emptyset$, we have 
    \begin{align}
\#S_{2,1/\alpha}=\sum\limits_{n=1}^{k/2}\left\lfloor{\dfrac{4n(k-m)}{2k+1}}\right\rfloor-\sum\limits_{n=1}^{k/2}\left\lfloor{\dfrac{(4n-2)(k-m)}{2k+1}}\right\rfloor\nonumber.
    \end{align}
Similarly, we can show that
    \begin{align}
\#S_{1,1/\alpha}=\sum\limits_{n=1}^{k/2}\left\lfloor{\dfrac{(4n+1)(k-m)}{2k+1}}\right\rfloor-\sum\limits_{n=1}^{k/2}\left\lfloor{\dfrac{(4n-1)(k-m)}{2k+1}}\right\rfloor\nonumber
    \end{align}
using the fact that $I_{11}=\emptyset$. The above formulas are obtained by counting the number of lattice points inside and on the triangle. To determine the number of lattice points within and on the triangle, we refer to (\cite{analytic}, p- 33). Since $S_{1,1/\alpha}$ and $S_{2,1/\alpha}$  are disjoint, we have 
    \begin{align}\label{cardS}
\#S_{1/\alpha}=&\sum\limits_{n=1}^{k/2}\left\lfloor{\tfrac{4n(k-m)}{2k+1}}\right\rfloor-\sum\limits_{n=1}^{k/2}\left\lfloor{\tfrac{(4n-2)(k-m)}{2k+1}}\right\rfloor+\sum\limits_{n=1}^{k/2}\left\lfloor{\tfrac{(4n+1)(k-m)}{2k+1}}\right\rfloor-\sum\limits_{n=1}^{k/2}\left\lfloor{\tfrac{(4n-1)(k-m)}{2k+1}}\right\rfloor\nonumber\\
=&\sum\limits_{n=1}^{2k+1}\varphi(n)\left\lfloor{\dfrac{n(k-m)}{2k+1}}\right\rfloor=\sum\limits_{n=1}^{2k}\varphi(n)\left\lfloor{\dfrac{n(k-m)}{2k+1}}\right\rfloor+k-m,
    \end{align}
    where $\varphi:J_{2k+1}\rightarrow \{-1,1\}$ is a function defined by
    \begin{align}\label{varphi}
        \varphi(n)=\begin{cases}
            1,&\text{if }n\equiv 0\text{ or }1\pmod{4},\\
            -1,&\text{if }n\equiv 2\text{ or }3\pmod{4}.
        \end{cases}
    \end{align}
\begin{figure}
\centering
\includegraphics[width=10.5cm]{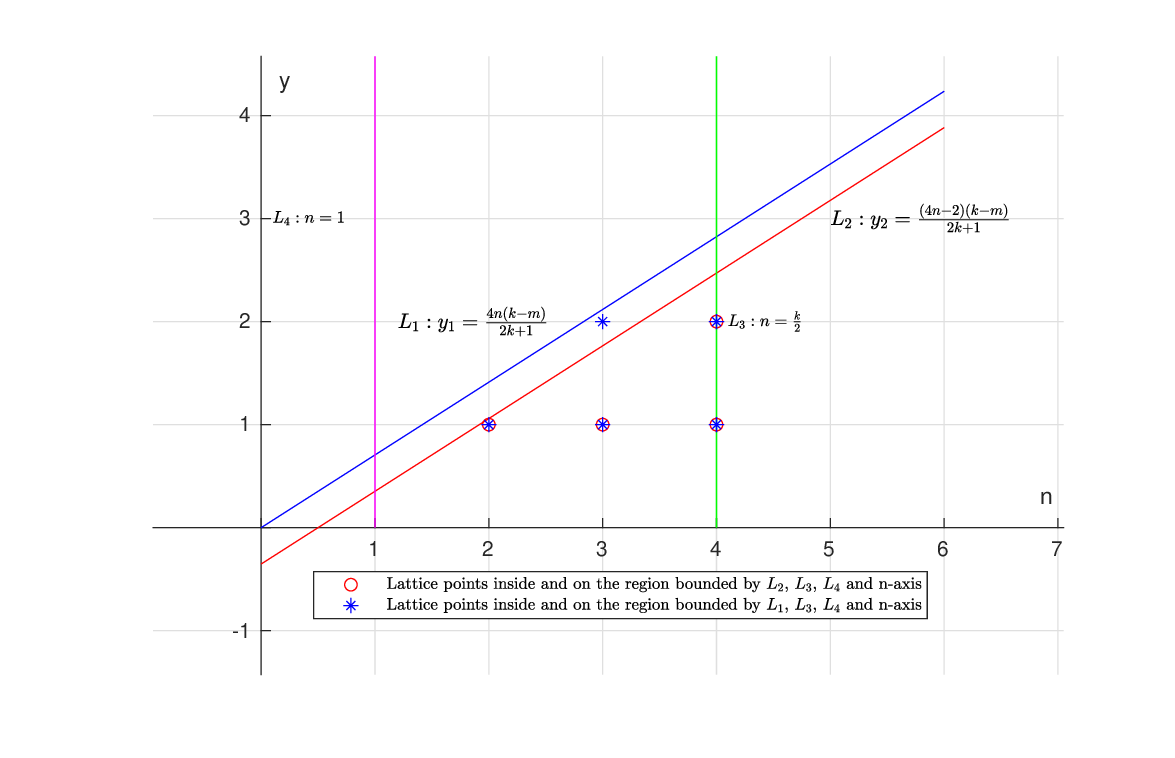}
\caption{The number of lattice points inside and on the region bounded by $L_1$, $L_2$, $L_3$, and $L_4$. The graph is actually drawn by choosing $m=5$ and $k=8$.}
\label{cardfig}
\end{figure}
    Since $k$ is even, $\varphi(2k+1-n)=\varphi(n),~\text{for all}~n\in J_{2k+1}.$
    After replacing $n$ by $2k+1-n$ and using $\floor{-x}=-1-\floor{x}$ for $x\notin \mathbb{Z}$, \eqref{cardS} becomes
    \begin{align*}     \#S_{1/\alpha}=~&\sum\limits_{n=1}^{2k}\varphi(2k+1-n)\left\lfloor{\dfrac{(2k+1-n)(k-m)}{2k+1}}\right\rfloor+k-m\\
=~&\sum\limits_{n=1}^{2k}\varphi(n)\left\lfloor{k-m-\dfrac{n(k-m)}{2k+1}}\right\rfloor+k-m\\
=~&\sum\limits_{n=1}^{2k}\varphi(n)\left(k-m-1-\left\lfloor{\dfrac{n(k-m)}{2k+1}}\right\rfloor\right)+k-m\\
=~&(k-m-1)\sum\limits_{n=1}^{2k}\varphi(n)-\sum\limits_{n=1}^{2k+1}\varphi(n)\left\lfloor{\dfrac{n(k-m)}{2k+1}}\right\rfloor+2(k-m)\\
=~&-\#S_{1/\alpha}+2(k-m),
\end{align*}
which implies that $\#S_{1/\alpha}=k-m$. 

$(ii)$ If $2m+1<\dfrac{1}{a}\le 2m+1+\dfrac{k-m}{2k+1},$
arguing as in $(i)$, we can show that $\#T_{1/\beta}=k-m$, when $\beta=2m+1+\dfrac{k-m}{2k+1}$. Therefore, $\#T_a\le k-m$. 
\end{pf}

\begin{lem} \label{mthm+1throws}
 $A_{mn}=A_{m+1,n}$, for $n=1,2,\dots,k$.
\end{lem}
\begin{pf}
Notice that
\begin{align*}
 X_{mn}=\dfrac{m(2k+1)-2n(2m+1)}{2(2m+1)(2k+1)}= \begin{cases}
        +ve,\text{ if } n=1,2,\dots,k/2-1,\\
        -ve,\text{ if } n=k/2,
    \end{cases}
\end{align*}
 and
 \begin{align*}
    X_{m,2k+1-n}=\dfrac{m(2k+1)-2(2k+1-n)(2m+1)}{2(2m+1)(2k+1)}
    \le~&\dfrac{m(2k+1)-2(k+1)(2m+1)}{8km+4k+4m+2}\\
    =~&-\dfrac{2km+2k+3m+2}{8km+4k+4m+2}<0,
 \end{align*}
 for all $n=k/2+1,\dots,k$. Further, \eqref{rangeofa}
 yields
 $$2a(2m+1)X_{m,2k+1-2n}\ge \dfrac{-2km-2k-3m-2}{4mk+k+3m+1}>-1.$$
Therefore, we obtain from Lemmas \ref{lem2} and \ref{lem3} that
\begin{align}\label{Amn}
   A_{mn}=\begin{cases}      
2n/b,&\text{if }n=1,\dots,k/2-1,\\
2am,&\text{if }n=k/2,\\
(2k+1-2n)/b,&\text{if }n=k/2+1,\dots,k.
\end{cases} 
\end{align}
Notice that
$X_{m+1,2k+1-n}$
\begin{align*}
    \hspace{1.5cm}=\dfrac{(m + 1)(2k + 1) - 2(2k+1-n)(2m + 1)}{2(2m+1)(2k+1)}\le-\dfrac{4km+k+7m+3}{8km+4k+4m+2}<0,
\end{align*}
 for all $n=1,2,\dots,k/2$. When  $n=1,2,\dots,k/2-1$, we have
\begin{align*}
  0>2a(2m+1)X_{m+1,2k+1-n}=~&\dfrac{a}{2k+1}\big[(m + 1)(2k + 1) - 2(2k+1-n)(2m + 1)\big]\\
  \le& -\dfrac{4km+k+7m+3}{4mk+k+3m+1}<-1
\end{align*}
and 
\begin{align*}
    0<2a(2m+1)(1+X_{m+1,2k+1-n})=~&\dfrac{a}{2k+1}\big[(m + 1)(2k + 1)+2n(2m + 1)\big]\\
    \le~&\dfrac{4km+3k-3m-1}{4km+k+3m+1}<1,
\end{align*}
because $4km+k+3m+1-(4km+3k-3m-1)=6m-2k+2>0$. By Lemma \ref{lem4},
\begin{align}\label{m+1firsthalf}
  A_{m+1,n}=2n/b,~n=1,\dots,k/2-1. 
\end{align}
When $2m+1-\dfrac{k-m}{2k+1}\le \dfrac{1}{a}\le 2m+1$ and $n=k/2$, we obtain from \eqref{rangeofa} that
\begin{align*}
0>2a(2m+1)X_{m+1,3k/2+1}=~&\dfrac{a}{2k+1}\big[(m+1)(2k+1)-2(3k/2+1)(2m+1)\big]\\
    =~&-\dfrac{a}{2k+1}(4mk+k+3m+1)\ge -1
\end{align*}
and 
\begin{align*}
2a(2m+1)\left(1+X_{m+1,3k/2+1}\right)=~&\dfrac{a}{2k+1}(4mk+3k+m+1)\\
 \ge~& \dfrac{4mk+3k+m+1}{4mk+2k+2m+1}>1.
\end{align*}
If $-Y<X_{m+1,3k/2+1}<0$, then
$$A_{m+1,k/2}=2ma$$
and if $X_{m+1,3k/2+1}=-Y,$ \textit{i.e.}, $\dfrac{1}{a}=2m+1-\dfrac{k-m}{2k+1}$, then
$$A_{m+1,k/2}=1-a(m+1)+\dfrac{k}{2b}=\dfrac{4km+2m}{4mk+3k+m+1}=2ma.$$
Similarly, we can show $A_{m+1,k/2}=2ma$, when $2m+1<\dfrac{1}{a}\le 2m+1+\dfrac{k-m}{2k+1}$.
Therefore,
\begin{align}\label{m+1k/2}
    A_{m+1,k/2}=2ma.
\end{align}
For $n=k/2+1,\dots,k$, we have
\begin{align*}
    X_{2m+1-s,n}=~&\dfrac{(m+1)(2k+1)-2(2k+1-n)(2m+1)}{2(2m+1)(2k+1)}\\
    =~&-\dfrac{6km+2k-m-1}{2(2m+1)(2k+1)}<0.
 \end{align*}
Therefore,
\begin{align}\label{m+1lasthalf}
   A_{m+1,n}=
    (2k+1-2n)/b,~n=k/2+1,\dots,k,
\end{align}
from Remark \ref{remlem5}.
Altogether equalities \eqref{Amn}, \eqref{m+1firsthalf}, \eqref{m+1k/2}, and \eqref{m+1lasthalf} prove our result.
\end{pf}
\begin{lem}\label{rank} 
Let $A$ be  a  $2m\times k$ matrix whose $(s,n)$-th entry is 
$$A_{sn}=\Phi_{sn}(0,0)-\Phi_{s,2k+1-n}(0,0),~s=1,2,\dots,2m, ~n=1,2,\dots,k.$$  Then
$\rank(A)\le k-1.$
\end{lem}
\begin{pf}
Let
$$\textbf{S}=\big\{s\in S_a\cup T_a: A_{s,n}\ne A_{2m+1-s,n},\text{for some } n=1,2,\dots,k\big\}.$$
Then either $\textbf{S}\subseteq S_a$ or $\textbf{S}\subseteq T_a$.\\

\noindent
\underline{\textbf{Case: 1}} When $\textbf{S}\subseteq S_a$. For $n=k/2$, $m$ satisfies the right inequality (in fact equality) that appears in \eqref{gauss2}. Hence $m\in S_{1/\alpha}$, where $\alpha=2m+1-\tfrac{k-m}{2k+1}$.
Consequently, we obtain from  Lemmas \ref{cardslem} $(i)$ and \ref{mthm+1throws} that $\#\textbf{S}\le k-m-1$. 
So, there are at least $2m-k+1$ pairs of identical rows in $A$ by Lemmas \ref{lem2} - \ref{lem5}. Hence $\rank(A)\le 2(k-m-1)+2m-k+1=k-1.$\\

\noindent
\underline{\textbf{Case: 2}} When $\textbf{S}\subseteq T_a$. Using a similar argument as in Case $1$, we can conclude $\rank(A)\le k-1$ from Lemmas \ref{cardslem} $(ii)$ and \ref{mthm+1throws}. 
\end{pf}

To prove the Conjecture \ref{conjecture2}, let us assume that
$ p = 2k + 1$ and $q =4m$. The condition that $\gcd(p,q)=1$ is necessary in Conjecture \ref{conjecture2}.
For example, if we choose $m=3$ and $k=4$, then $ab=\tfrac{9}{12}$, $\gcd(9,12)=3$ but $\mathcal{G}(Q_2,1/2,3/2)$ is a frame (see \cite{ofsb, ourgabor}). Let us define
$$\widetilde{X}_{sn}=\dfrac{s}{4m}-\dfrac{n}{2k+1},\text{ for }s = 0, \dots, 4m - 1\text{ and }n = 1, \dots , 2k.$$ 
We can easily verify that $-1<\widetilde{X}_{sn}<1$ and $\widetilde{X}_{sn}\ne 0$ for all $s = 1, \dots , 4m -1$ and $n = 1, \dots, 2k$. To distinguish the symbol matrix from $\Phi(0,0)$, we denote it here by $\widetilde{\Phi}(0,0)$.
Now arguing as in Lemma \ref{lem1}, we can show that
\begin{align}\label{newphi00}
   \widetilde{\Phi}_{sn}(0,0)=\begin{cases}
\sum\limits_{l=-1}^0 Q_2\big(2(l+\widetilde{X}_{sn})\big), & \text{if}~ 0<\widetilde{X}_{sn}<1,\\
\sum\limits_{l=0}^1 Q_2\big(2(l+\widetilde{X}_{sn})\big), & \text{if} -1<\widetilde{X}_{sn}<0,
\end{cases} 
\end{align}
for all $s=1,\dots,4m-1$ and $n=1,\dots,2k$. 

Let us define
\begin{eqnarray}
\widetilde{A}_{sn}=\widetilde{\Phi}_{sn}(0,0)-\widetilde{\Phi}_{s,2k+1-n}(0,0),~\text{for}~s = 0,\dots,4m-1, ~n = 1,\dots, k.
\end{eqnarray}
Using a similar argument as in Lemmas \ref{02m+1th0} - \ref{lem5}, we can prove the following Lemmas.
\begin{lem}\label{02m+1th0new}
    \noindent
    \begin{enumerate}
        \item [$(i)$] $\widetilde{A}_{0n}=\widetilde{A}_{2m,n}=0$, for $n=1,2,\dots,k$.
        \item[$(ii)$] $\widetilde{A}_{s,n}=-\widetilde{A}_{4m-s,n}$,
        for $s=1,2,\dots,2m-1$ and  $n=1,2,\dots,k$.   
    \end{enumerate}    
\end{lem}
\begin{lem}\label{newlem}
    Assume that $k$ is even. Let $s = 1,\dots, m-1$. Then
    \begin{enumerate}
        \item [$(i)$] For $n = 1, 2, \dots, k/2$, we have 
$$\widetilde{A}_{sn}=\begin{cases}
{2n}/{b}, & \text{if}~ 0<\widetilde{X}_{sn}<\tfrac{1}{2},\\
s/m, & \text{if} -\tfrac{1}{2}<\widetilde{X}_{sn}<0.
\end{cases}$$
\item[$(ii)$] For 
$n = k/2+1, \dots, k$, we have
$$\widetilde{A}_{sn}=
\begin{cases}
{(2k+1-2n)}/{b}, & \text{if}~ -\tfrac{1}{2}<\widetilde{X}_{s,2k+1-n}<0,\\
s/m, & \text{if} -1<\widetilde{X}_{s,2k+1-n}<-\tfrac{1}{2}.
\end{cases} $$
\item[$(iii)$] For $n = 1, \dots, k/2$, we have 
$$
\widetilde{A}_{2m-s,n}=\begin{cases}
s/m, & \text{if}~ -\tfrac{1}{2}<\widetilde{X}_{2m-s,2k+1-n}<0,\\
{2n}/{b}, & \text{if} -1<\widetilde{X}_{2m-s,2k+1-n}<-\tfrac{1}{2}.
\end{cases}    
$$
\item[$(iv)$] For 
$n = k/2+1, \dots, k$, we have $$\widetilde{A}_{2m-s,n}=\begin{cases}
s/m, & \text{if}~ 0<\widetilde{X}_{2m-s,n}<\tfrac{1}{2},\\
{(2k+1-2n)}/{b}, & \text{if} -\tfrac{1}{2}<\widetilde{X}_{2m-s,n}<0.
\end{cases}$$
\end{enumerate}
\end{lem}

\begin{lem}\label{newrankA}
Assume that $k$ is even. Let $\widetilde{A}$ be  a  $(2m-1)\times k$ matrix whose $(s,n)$-th entry is 
$$\widetilde{A}_{sn}=\widetilde{\Phi}_{sn}(0,0)-\widetilde{\Phi}_{s,2k+1-n}(0,0),~s=1,\dots,2m-1, ~n=1,\dots,k.$$  Then
$\rank(\widetilde{A})\le k-1.$
\end{lem}
\begin{pf}
    Notice that
    \begin{align*}
        -1<2\widetilde{X}_{2m-s,2k+1-n}<0\Longleftrightarrow& -1<2\left(\dfrac{2m-s}{4m}-\dfrac{2k+1-n}{2k+1}\right)<0\\
        \Longleftrightarrow&-1<2\left(\dfrac{1}{2}-\dfrac{s}{4m}-1+\dfrac{n}{2k+1}\right)<0\\
        \Longleftrightarrow&1>2\left(\dfrac{1}{2}+\dfrac{s}{4m}-\dfrac{n}{2k+1}\right)>0\\
        \Longleftrightarrow &-1<2\widetilde{X}_{sn}<0
    \end{align*}
    and 
    \begin{align*}
        -2<2\widetilde{X}_{2m-s,2k+1-n}<-1\Longleftrightarrow& -2<2\left(\dfrac{2m-s}{4m}-\dfrac{2k+1-n}{2k+1}\right)<-1\\
        \Longleftrightarrow&-2<2\left(\dfrac{1}{2}-\dfrac{s}{4m}-1+\dfrac{n}{2k+1}\right)<-1\\
        \Longleftrightarrow&2>2\left(\dfrac{1}{2}+\dfrac{s}{4m}-\dfrac{n}{2k+1}\right)>1\\
        \Longleftrightarrow &0<2\widetilde{X}_{sn}<1.
    \end{align*}
    Therefore, from Lemma \ref{newlem} $(i)$ and $(iii)$, $$\widetilde{A}_{sn}=\widetilde{A}_{2m-s,n},\text{ for }s=1,\dots,m-1,~n=1,\dots,k/2.$$ 
    Similarly, from Lemma \ref{newlem} $(ii)$ and $(iv)$ we can show that $$\widetilde{A}_{sn}=\widetilde{A}_{2m-s,n},\text{ for }s=1,\dots,m-1,~n=k/2+1,\dots,k.$$ Therefore, $\widetilde{A}_{sn}=\widetilde{A}_{2m-s,n}$ for all $s=1,\dots,m-1$ and $n=1,\dots,k$. Hence 
    $$\rank{\widetilde{A}}\le m-1+1=m\le k-1.$$
\end{pf}

\begin{rem}
  We can prove that Lemmas $\ref{lem2}-\ref{newrankA}$ are also valid for odd $k$ by separating the range set of $n$ into two sets $\{1,\dots,(k-1)/2\}$ and $\{(k+1)/2,\dots,k\}$ and suitably modifying $S_{ja}$ sets. We also verified them.
\end{rem}

\section{Proof of Theorem \ref{main thm}}
In Section 3, we built up essential tools to prove our main results. Now we
are ready to prove Theorem \ref{main thm} in this section.
\newline
\underline{\textbf{Proof of Conjecture \ref{conjecture}: }}
For simplicity, we write $\Phi_{sn}(0,0)\equiv \Phi_{sn}.$
Let $C_l$ be the $l$-th column of $\Phi(0,0)$. Now replacing $C_l$ by $C_l-C_{2k+1-l}$, for $l=1,\dots,k$, by Lemmas \ref{02m+1th0} and \ref{reflection}, $\Phi(0,0)$ is column equivalent to the following matrix:
\begin{equation}\label{clumnequivphi}
   \left[\setlength{\arraycolsep}{0.29cm}\begin{array}{c|c|c}
     \Phi_{00} &\begin{matrix}
         0&0&\dots&0
     \end{matrix}&\begin{matrix}
     \Phi_{0,k+1}&\dots&\Phi_{0,2k}\end{matrix}\\
     \hline
\begin{matrix}
    \Phi_{10}\\
    \vdots\\
    \Phi_{2m,0}
\end{matrix}&A&\begin{matrix}
    \Phi_{1,k+1}&\dots&\Phi_{1,2k}\\
    \vdots&\ddots&\vdots\\
    \Phi_{2m,k+1}&\dots&\Phi_{2m,2k}
\end{matrix}\\
\hline
\Phi_{2m+1,0} &\begin{matrix}
         0&0&\dots&0
     \end{matrix}&\begin{matrix}
     \Phi_{2m+1,k+1}&\dots&\Phi_{2m+1,2k}\end{matrix}\\
     \hline
     \begin{matrix}
    \Phi_{2m+2,0}\\
    \vdots\\
    \Phi_{4m+1,0}
\end{matrix}&-B&\begin{matrix}
    \Phi_{2m+2,k+1}&\dots&\Phi_{2m+2,2k}\\
    \vdots&\ddots&\vdots\\
    \Phi_{4m+1,k+1}&\dots&\Phi_{4m+1,2k}
\end{matrix}
    \end{array}
    \right],
\end{equation}
where $A$ is the $2m\times k$ matrix defined in Lemma \ref{rank} and $B$ is the reflection matrix of $A$ which is defined as $B_{sn}=A_{2m+1-s,n}$, for $s=1,2,\dots ,2m$ and $n=1,2,\dots,k$.

Since $\rank(A)\le k-1$, there exist scalars $\gamma_1, \gamma_2,\dots,\gamma_k$, not all zero such that  
$\gamma_1 A_1+\cdots+\gamma_k A_k=0$, where $A_1$, $A_2,\dots,A_k$ are columns of the matrix $A$. Without loss of generality, we assume that $\gamma_k\ne 0$. If we replace the $k-$th column of matrix $A$ by $A_k+\tfrac{1}{\gamma_k}(\gamma_1 A_1+\cdots+\gamma_{k-1}A_{k-1})$, then $A$ is column equivalent to a matrix whose $k-$th column is the zero vector.
Since $B$ is the reflection matrix of $A$, after applying the same column operation on $B$, $B$ is column equivalent to a matrix whose $k-$th column is the zero vector. From these observations, if we replace the $k-$th column of matrix in \eqref{clumnequivphi} by $C_k+\tfrac{1}{\gamma_k}(\gamma_1 C_1+\cdots+\gamma_{k-1}C_{k-1})$, then the matrix in \eqref{clumnequivphi} is column equivalent to a matrix whose $k$-th column is a zero vector. Therefore, $\Phi(0,0)$ has linearly dependent columns and hence $\Phi^{*}\Phi$ is not invertible.  Hence our proof follows from Theorem \ref{framewithrank}.
\newline
\underline{\textbf{Proof of Conjecture \ref{conjecture2}: }} The proof of Conjecture \ref{conjecture2} is carried out in the same way as in the proof of Conjecture \ref{conjecture}.
The main difference in this case is that the $4m\times (2k+1)$ symbol matrix $\widetilde{\Phi}(0,0)$ is column equivalent to a matrix which is similar to that in \eqref{clumnequivphi} with proper modifications on entries; $\widetilde{A}$  will take the place of $A$ which is defined in Lemma \ref{newrankA}, and $B$ will be replaced by $\widetilde{B}=\big[\widetilde{A}_{2m-s,n}\big]_{s=1,n=1}^{2m-1,k}$.
Now arguing as in the proof of the Conjecture \ref{conjecture}, we can conclude our result from Lemma \ref{newrankA}.

\end{document}